\documentclass[a4paper,11pt]{article}
\usepackage{amsmath}
\usepackage{amssymb}
\usepackage {mathrsfs}
\usepackage{graphics}
\usepackage{amscd}
\begin{document}
\newcommand{\projsp}[1]{\mathbb{P}^#1}
\newcommand{\qed}{$\blacksquare$}
\newtheorem{lemma}{Lemma}
\newtheorem{theorem}{Theorem}
\newtheorem{remark}{Remark}
\newtheorem{corollary}{Corollary}
\newtheorem{example}{Example}
\title{ Double Rational Normal Curves with Linear Syzygies}
\author{Nicolae Manolache \\ Institute of Mathematics of the Romanian Academy, \\
P.O. Box 1--764, RO--70700 Bucharest, Romania, \\
email: Nicolae.Manolache@imar.ro\\
current address: Fachbereich 6 Mathematik, \\ Universit\"at Oldenburg, Pf 2503,
D-26111 Oldenburg, Germany, \\
email: nicolae.manolache@uni-oldenburg.de}
\date{}
\maketitle
\begin{abstract}
In this note we are looking after nilpotent projective curves without embedded 
points, which have rational normal  curves of degree $d$ as support, 
are defined (scheme-theoretically) by quadratic equations, have degree $2d$ and 
have only linear syzygies. We show that, as expected, no such curve does exist in 
$\projsp{d}$, and then consider doublings in a bigger ambient space. The simplest
and trivial example is that of a double line in the plane. We show 
that the only possibility is to take a certain doubling in the sense of
Ferrand (cf \cite{RefF}) of rational normal  curves in $\projsp{d}$ 
embedded further linearly in $\projsp{{2d}}$.
These double curves have the 
Hilbert polynomial $H(t)=2dt+1$, i.e. they are in the Hilbert scheme of the 
rational normal curves of degree $2d$. Thus, it turns out that they are natural
generalizations of the double line in the plane as a degenerated conic.
The simplest nontrivial
example is the curve of degree $4$ in $\projsp{4}$, defined by the ideal
$(xz-y^2,\ xu-yv,\ yu-zv,\ u^2,\ uv,\ v^2)$. The double rational normal curve allow the
formulation of a {\em Strong Castelnuovo Lemma} in the sense of \cite{RefG},
for sets of points and double points. In the last section we mention some plethysm formulae
for symmetric powers.
\end{abstract} 
\section{Introduction}
\label{sec:0}
The syzygies for Veronese or Segre embeddings are not known, with some simple 
exceptions.
In very few cases they are "pure" (cf. \cite{RefBM1} for a list, in which  the simple case
of the Segre embedding  $\projsp{1} \times \projsp{2} \rightarrow \projsp{5}$ is missing,
or use \cite{RefS} to produce the same list, as "extremal rings of format 2")
and no other cases are given  explicitely in the
literature. Already beginning with the Veronese embeddings $v_d(\projsp{2})$ the
syzygies are not linear, excepting some small $d$'s. 
It was shown in \cite{RefB} that
the first $3d-3$ syzygies are linear and in \cite{RefOP} that exactly at the step $3d-2$
(if $d\ge 3$) they  are failing to be linear. The case $d=3$ is easy and already known
(cf \cite{RefBM1}), the resolution being pure. The simplest case when the syzygies are
not pure is the Veronese embedding $v_4$ of $\projsp{2}$ into
$\projsp{{14}}$, given by monomials of degree $4$ (cf. {\bf Lemma} \ref{impure}). 
The syzygies of the rational normal  line ($\projsp{1}$ embedded via
the Veronese morphism $v_d$ in  $\projsp{d}$) are linear and are well
understood, being given by the Eagon-Northcot complex.
The nilpotent curves which have a rational normal  curve as support
are far from having quadratic equations or linear syzygies, even when we consider
"mild" nilpotency (locally Cohen-Macaulay or even complete intersections).
In this paper we construct, for each d, a family of double structures on the normal rational 
curve of degree $d$ embedded in $\projsp{{2d}}$, which have linear
syzygies. More precisely, we consider the embedding of $\projsp{1} \to \projsp{{2d}}$ 
given by  $v_d$ composed with a linear embedding 
$\projsp{d} \hookrightarrow \projsp{{2d}}$
and double conveniently this curve by the Ferrand's method (cf. \cite{RefF}).

In the last section, which strictly speaking does not interact with the rest of this note, 
one gives some "plethysm formulae" for symmetric powers, as consequences of minimal 
resolutions of Veronese embeddings.

\section{Preliminaries}
\label{sec:1}
We fix an algebraically closed field $k$ and all our schemes are algebraic schemes 
over it. 

Let $X$ be a subscheme in a projective space $\projsp{n}$ and $\mathscr{I}_X$
its sheafideal in $\mathcal{O} =\mathcal{O}_{\projsp{n}}$.
By definition, the graded ideal of $X$ is $I_X=\oplus _{\ell \ge 0}H^0(\mathscr{I} 
(\ell))$. To a minimal resolution of the graded ideal $I:=I_X$, over the graded ring 
of polynomials $R=k[X_0,\ldots ,X_n]$, let say of the form:
\[
0 \leftarrow I \leftarrow \oplus_{j=1}^{j=r_1} R(-d_{1j})^{b_{1j}} \leftarrow 
\oplus_{j=1}^{i=r_2} R(-d_{2j})^{b_{2j}}\leftarrow \ldots  
\]
written also:
\[
0 \leftarrow I \leftarrow \oplus_{j=1}^{j=r_1}b_{1j} R(-d_{1j})\leftarrow 
\oplus_{j=1}^{j=r_2}b_{2j}R(-d_{2j})\leftarrow \ldots  
\]
one associates canonically a resolution of the sheaf of ideals $\mathscr{I}:=\mathscr{I}_X$
of the type:
\[
0 \leftarrow \mathscr{I} \leftarrow \oplus_{j=1}^{j=r_1}b_{1j} \mathcal{O} (-d_{1j})
\leftarrow \oplus_{j=1}^{j=r_2}b_{2j}\mathcal{O} (-d_{2j})\leftarrow \ldots  
\quad \quad ,
\]
which we shall call also \emph{ minimal resolution of } $\mathscr{I}$. 

The corresponding resolution of $\mathcal{O} _X$ :
\[
0 \leftarrow \mathcal{O} _X \leftarrow \mathcal{O} \leftarrow 
\oplus_{j=1}^{j=r_1}b_{1j} \mathcal{O} (-d_{1j})\leftarrow 
\oplus_{j=1}^{j=r_2}b_{2j}\mathcal{O} (-d_{2j})\leftarrow \ldots  
\]
will be called  \emph{ a minimal resolution of } $\mathcal{O} _X$, although it does not 
corresponds in general to a minimal resolution of the associated graded $R$--module
$\oplus H^0(\mathcal{O} (\ell))$.
The $\mathcal{O}$--homomorphisms in these resolutions correspond to matrices
of homogeneous polynomials in $X_i$, whose degrees are determined by the numbers 
$d_{ij}$.

To a (minimal) resolution of a graded $R$--module $M$, let say of the form:
\[
0 \leftarrow M \leftarrow \oplus_{j=1}^{j=r_1}b_{0j} R(-d_{0j})\leftarrow 
\oplus_{j=1}^{j=r_2}b_{1j}R(-d_{1j})\leftarrow \ldots  
\]
one associates the polynomial in two variables:
\[
P_M(x,t)=\sum_{ij} b_{ij}t^{d_{ij}}x^i
\]
which depends only on $M$ (cf. \cite{RefG}) and contains the whole information as far as 
the ``Betti numbers'' $b_{ij}$ and the ``shiftings'' of degrees $d_{ij}$ are 
concerned.

We shall use also a notation  for a minimal resolution considered in \cite{RefG}: 
\[
0\leftarrow M \leftarrow \bigoplus _j {\bf M}_{0j}(M) \otimes R(-j) \leftarrow \ldots
\leftarrow \bigoplus _j {\bf M} _{ij}(M) \otimes R(-j) \leftarrow 
\]
where ${\bf M}_{ij}(M)$ are vector spaces, only finiteley many different from $0$.

The following lemma should be very well known, but we know no refference for it.
\begin{lemma} Let $X$ be a subscheme in $\projsp{n}$ and consider $\projsp{n}$
embedded linearly in $\projsp{{n+m}}$. If
$P_{X,n}(x,t) = $
{\rm the polynomial associated to the minimal resolution 
of } $\mathcal{O} _X$ {\rm  in } $\projsp{n}$,
then the polynomial
associated to the minimal resolution in $\projsp{{n+m}}$ is
\[
P_{X,n+m}(x,t)= P_{X,n}(x,t)(1+xt)^m
\]
\end{lemma}
Proof. It is enough to prove this lemma for $m=1$. 
Observe that a minimal resolution of the graded module $M$ over $R=k[X_0,\ldots ,
X_n]$ :
\[
0\leftarrow M \leftarrow L_0 \xleftarrow{A}L_1 \xleftarrow{B} L_2 
\xleftarrow{C} L_3 \leftarrow \ldots \quad ,
\]
where $L_i=\oplus_j R(-d_{ij})$,
gives the following minimal resolution of $M$ as $R' = k[X_0, \ldots , X_{n+1}]=
R[X_{n+1}]$--module:
\begin{eqnarray}
0\leftarrow M \leftarrow L_0' \xleftarrow{\begin{pmatrix} A & -X 
\end{pmatrix}} L_1'\oplus L_0'(-1) 
\xleftarrow{\begin{pmatrix} B & X \\ 0 & A 
\end{pmatrix}} \nonumber \\ 
\leftarrow L_2'\oplus L_1'(-1) 
\xleftarrow{\begin{pmatrix}C & -X\\ 0 & B\end{pmatrix}} L_3'\oplus L_2'(-1) 
\leftarrow \ldots \quad ,\nonumber
\end{eqnarray}
where $L_i'= \oplus _j R'(-d_{ij})$ and $X=\hbox{\rm multiplication by }X_{n+1}$\ .
\qed

As it is very well known, the syzygies of the rational normal  curve are given by
an Eagon-Northcott complex. Let $X$ be the rational normal  curve of degree $d$
in $\projsp{d}$, i.e. $X= {\rm \hbox{ the image of }}\projsp{1}  {\rm \hbox{ under 
the Veronese embedding } } v_d$. 
Then one has the following minimal resolution (cf. e.g. \cite{RefBM1}, \cite{RefG}):
\begin{eqnarray}
0 \leftarrow {\mathcal{O} }_X \leftarrow {\mathcal{O} } \leftarrow {d \choose 2}
{\mathcal{O}}(-2)\leftarrow 2{d \choose 3}{\mathcal{O}}(-3)\leftarrow 3{d \choose 4}
{\mathcal{O}}(-4)\leftarrow & \nonumber \\
\ldots \leftarrow \ell {d \choose {\ell +1}}{\mathcal{O}}(-\ell -1)\leftarrow \ldots 
\leftarrow (d-1){d \choose d}{\mathcal{O}}(-d)\leftarrow 0 & \nonumber
\end{eqnarray}

For the following considerations  we shall need also the resolutions of all 
line bundels on $X$. As $X\cong \projsp{1}$, $Pic(X)\cong \mathbb{Z}$. We denote by
$\mathcal{T}$ a line bundle on $X$ corresponding to the generator 
$\mathcal{O} _{\projsp{1}}(1)$
of $Pic(X)$. $\mathcal{T}$ is not induced from any line bundle on $\projsp{d}$, 
but all powers
$\mathcal{T}^{\ell d}$ are, namely $\mathcal{T}^{\ell d}=v_d^*(O_{\projsp{d}}(\ell))$. 
For us it is sufficient to know the minimal resolution of $\mathcal{T}$, $\mathcal{T}^2$, 
\ldots , $\mathcal{T}^{d-1}$, but the general case is equally easy.
\begin{lemma}
If $\mathcal{T}$ is a line bundle on the rational normal curve $X$ in $\projsp{d}$, which 
corresponds to $\mathcal{O}_{\projsp{1}}(1)$, then the minimal resolutions of
$\mathcal{T}^j$ in $\projsp{d}$ ($j=1,\ldots d-2$), are of the shape:
\begin{eqnarray}
0\leftarrow \mathcal{T}^j  \leftarrow (j+1)\mathcal{O} \leftarrow j{d \choose 1}
\mathcal{O}(-1) 
\leftarrow (j-1){d \choose 2}\mathcal{O}(-2)\leftarrow \nonumber \\ 
\leftarrow (j-2){d \choose 3}\mathcal{O}(-3) \leftarrow \ldots  
\leftarrow {d \choose j} \mathcal{O}(-j) \leftarrow {d \choose j+2}\mathcal{O}(-j-2)
\leftarrow \nonumber \\
\leftarrow 2{d \choose j+3}\mathcal{O} (-j-3) \leftarrow 3{d \choose j+4}\mathcal{O} (-j-4)
\leftarrow \ldots \leftarrow \nonumber \\
\leftarrow (d-j-1){d \choose d}\mathcal{O}(-d) \leftarrow 0 \nonumber
\end{eqnarray}
and the resolution of $\mathcal{T} ^{d-1}$ is of the shape:
\begin{eqnarray}
0\leftarrow \mathcal{T}^{d-1}  \leftarrow d\mathcal{O} \leftarrow (d-1){d \choose 1}
\mathcal{O}(-1) 
\leftarrow (d-2){d \choose 2}\mathcal{O}(-2)\leftarrow \nonumber \\ 
\leftarrow (d-3){d \choose 3}\mathcal{O}(-3) \leftarrow \ldots  
\leftarrow {d \choose d-1} \mathcal{O}(-d+1) \leftarrow 0 \nonumber
\end{eqnarray}
\end{lemma}
Proof. Standard exercise. Observe that only for $j=d-1$ all the syzygies are linear. \qed

\section{Doubling rational normal  curves of degree $d$}
\label{sec:2}
We recall shortly the Ferrand's method of doubling (cf. \cite{RefF}, \cite{RefBF} or \cite{RefM1}).
Given a locally Cohen-Macaulay curve $X$ embedded in a regular scheme $\mathbb{P}$, 
let $N_X$ be the conormal sheaf of $X$ in $\mathbb{P}$, i.e. $N_X=\mathscr{I} /\mathscr{I} 
^2$, where $\mathscr{I}$ is the sheaf of ideals of $X$ in $\mathbb{P}$. Let $\mathcal{L}$ be 
an invertible bundle on $X$,  $\omega = \omega _X$ be the dualizing sheaf of $X$, 
and $p:N_X\to \omega  \otimes  \mathcal{L}$ a surjective homomorphism. 
Then the kernel of $p$ will have the form $\mathscr{J}/\mathscr{I}^2$, where $\mathscr{J}$ 
is an ideal in ${\mathcal{O}}_\mathbb{P}$, which defines a scheme $Y$ having the support 
$X$, but with double multiplicity in every point of $X$. One has the following exact 
sequence:
\[
0 \rightarrow \mathscr{J}/\mathscr{I}^2  \rightarrow N_X \rightarrow \omega \otimes 
 \mathcal{L} \rightarrow 0
\]
which gives the exact sequences:
\[
0 \rightarrow \mathscr{J} \rightarrow \mathscr{I} \rightarrow \omega \otimes 
 \mathcal{L} \rightarrow 0 \quad ,
\]
\[
0 \rightarrow \omega \otimes  \mathcal{L} \rightarrow \mathcal{O} _Y \rightarrow \mathcal{O} _X 
\rightarrow 0
\]
and 
\[
0 \rightarrow \omega \otimes \mathcal{L} \rightarrow \mathcal{O} _Y^*\rightarrow \mathcal{O} _X^* 
\rightarrow 0 
\quad .
\]
From the last exact sequence one obtains the exact sequence:
\[
H^1(\omega \otimes \mathcal{L}) \rightarrow Pic(Y) \rightarrow Pic(X) \rightarrow H^2(\omega 
\otimes  \mathcal{L})
\]
The curve $Y$ is locally Gorenstein and one shows easily that its dualizing sheaf is
\[ \omega _Y |_X\cong  \mathcal{L}^{-1} \quad .\]
In the case when the curve $X$ is Gorenstein, any locally free sheaf of rank one
can be written as $\omega \otimes \mathcal{L}$, with  $\mathcal{L}$ conveniently chosen, so 
that a doubling of $X$ is given by a surjection $ \mathscr{I} /\mathscr{I}^2 \to  \mathcal{L}$, 
$ \mathcal{L}$
a locally free sheaf of rank $1$. The above exact sequences are written then with 
$\mathcal{L}$ instead of $\omega \otimes  \mathcal{L}$.

We recall that
the normal bundle to the rational normal curve of degree $d\ge 3$, considered as a bundle 
on $\projsp{1}$,
is 
\[(d-1)\mathcal{O}_{\projsp{1}}(d+2)=\oplus _{j=1}^{d-1}\mathcal{O} _{\projsp{1}}(d+2)
\]
\begin{lemma}
There is no double structure on the rational normal  curve of degree d in $\projsp{d}$
whose homogeneous ideal is generated by quadrics and has only linear syzygies.
\end{lemma} 
Proof. Suppose that such a double curve does exist.  The case $d=2$ being trivial, suppose
$d\ge 3$. Let
\[
p:\mathscr{I}/\mathscr{I} ^2 \rightarrow  \mathcal{L}\quad, \quad  \mathcal{L}|_{\projsp{1}}=
\mathcal{O}_{\projsp{1}}(r), \quad r\ge -d-2
\]
be the data defining it.  Let $\mathscr{J}$ be the sheaf of ideals of $Y$ in $\projsp{d}$.
It is easy to see from the exact sequence:
\[
0\rightarrow J \rightarrow \oplus _{t\ge 0} H^0(\mathcal{L}(t)) \rightarrow 
\oplus _{t\ge 0}H^1(\mathcal{J}(t)) \rightarrow 0
\]
that the homogeneous ideal $J$ of $Y$ has $\mathrm{depth}(J)\ge 2$, 
hence the minimal resolution of $\mathscr{J}$ has length at most $d-1$:
\[
0\leftarrow \mathscr{J} \leftarrow b_1\mathcal{O} (-2) \leftarrow b_2\mathcal{O} (-3) \leftarrow
\ldots \leftarrow b_{d-1}\mathcal{O} (-d) \leftarrow b_d \mathcal{O} (-d-1) \leftarrow 0
\]
It follows:
\begin{eqnarray}
h^1(\mathscr{J}) = h^d(b_d\mathcal{O}(-d-1))=b_d  \nonumber \\
H^1(\mathscr{J}(1))\cong H^d(b_d\mathcal{O}(-d))=0 \nonumber
\end{eqnarray}
The exact sequence :
\begin{equation}
0\rightarrow \mathscr{J} \rightarrow \mathscr{I} \rightarrow  \mathcal{L} \rightarrow 0 
\label{0} \tag{$\ast$}
\end{equation}
gives:
\begin{align}
b_d = & h^1(\mathscr{J}) = h^0( L)  \tag{$\ast\ast$} \label{1}\\
0 = & h^1(\mathscr{J}(1))\cong H^0( L(1))=H^0(\mathcal{O} _{\projsp{1}}(d+r)) \quad .
\tag{$\ast\ast\ast$} \label{2}
\end{align}
(\ref{2})  shows that $r<-d$ and (\ref{1}) gives then $b_d=0$, i.e. $Y$ is
arithmetically Cohen-Macaulay in $\projsp{d}$. 
The resolution of $\mathscr{I}$ shows that $H^2(\mathscr{I} )=0$. 
Then the exact sequence (\ref{0})
gives $H^1( L)=0$, i.e. $r\ge -1$. So far we got $-1\le r<-d$, contradiction.
\qed

\section{Rational Normal Curves of Degree $d$ Embedded in $\projsp{{d+e}}$}
\label{sec:3}
Because of the result of the last section, it remains to look after doublings
of rational normal  curves lying in a linear subspace of the projective space.
To fix the notation, let $X$ be the image of the embedding:
\[
\projsp{1}\xrightarrow{Veronese}\projsp{d}\xrightarrow{linear}\projsp{{d+e}} \quad .
\]
$X$ is arithmetically Cohen-Macaulay in $\projsp{{d+e}}$ and, according to 
{\bf  Lemma 1}, the 
polynomial associated to the minimal resolution of §$X$ in $\projsp{{d+e}}$
is the polynomial corresponding to the embedding of $\projsp{1}$ in 
$\projsp{d}$ multiplied with $(1+tx)^e$.
\begin{theorem}
There exist double structures $Y \subset \projsp{{d+e}}$ on rational normal curves 
$X$ of degree $d$ in $\projsp{d} \subset \projsp{{d+e}}$ defined by quadratic equations and
having only linear syzygies, only when $e=d$ and the 
line bundle $ \mathcal{L}$ associated to the 
doubling is such that $ \mathcal{L}|_{\projsp{1}}=\mathcal{O} _{\projsp{1}}(-1)$.
The minimal resolution is then of the shape:
\begin{eqnarray}
0 \leftarrow {\mathcal{O} }_Y \leftarrow {\mathcal{O} } \leftarrow {2d \choose 2}
{\mathcal{O}}(-2)\leftarrow 2{2d \choose 3}{\mathcal{O}}(-3)\leftarrow 3{2d \choose 4}
{\mathcal{O}}(-4)\leftarrow & \nonumber \\
\ldots \leftarrow \ell {2d \choose {\ell +1}}{\mathcal{O}}(-\ell -1)\leftarrow \ldots 
\leftarrow (2d-1){2d \choose d}{\mathcal{O}}(-2d)\leftarrow 0 & \nonumber
\end{eqnarray}
i.e. is that of a rational normal curve of degree $2d$ in $\projsp{{2d}}$.
\end{theorem}
Proof. 
 Like in the proof of {\bf  Lemma 3}, one shows that a double structure $Y$
with linear syzygies should be arithmetically Cohen-Macaulay.
In this case finding the minimal
resolution of $Y$ (i.e. of $\mathscr{J}$) is equivalent to finding the minimal 
resolution of $\mathcal{O} _Y$, which should be of the shape:
\[
0\leftarrow \mathcal{O}_Y \leftarrow \mathcal{O} 
\leftarrow b_1 \mathcal{O} (-2) \leftarrow \ldots 
\leftarrow b_{d+e-1}\mathcal{O}(-d-e)\leftarrow 0
\]
It is easy to see, in fact like in the proof of {\bf Lemma 3}, that
the following vanishings take place:
\[
H^1(\mathscr{J} (\ell))=0, \quad H^2(\mathscr{J} (\ell))=0, 
\quad H^1( \mathcal{L} (\ell))=0, \quad (\ell \ge 0)\quad .
\]
It follows that, if $\mathcal{L} |_\projsp{1}=\mathcal{O}_{\projsp{1}}(r)$,
then $r\ge -1$ and the graded rings $R(X)$, $R(Y)$ associated to 
$\mathcal{O} _X$, $\mathcal{O} _Y$ and the graded module $M$ associated to $ \mathcal{L}$ fit in 
an exact sequence:
\[
0\leftarrow R(X) \leftarrow R(Y)  \leftarrow M \leftarrow 0 \quad ,
\]
which comes from the exact sequence:
\[
0\leftarrow \mathcal{O} _X \leftarrow \mathcal{O} _Y \leftarrow  \mathcal{L} \leftarrow 0
\]
Then we have for each $j$ the following exact sequence of vector spaces (cf. \cite{RefG}):
{\small
\[
0\leftarrow {\bf M}_{0j}(X) \leftarrow {\bf M}_{0j}(Y)\leftarrow {\bf M}_{0j}( \mathcal{L})
\leftarrow {\bf M}_{1j}(X) \leftarrow {\bf M}_{1j}(Y) \leftarrow  {\bf M}_{1j}( \mathcal{L})
\leftarrow \ldots ,
\]
}
where ${\bf M}_{0j}(X):={\bf M}_{0j}(R(X))$, ${\bf M}_{0j}(Y):={\bf M}_{0j}(R(Y))$.

In the following we shall denote by ${\bf m}_{ij}(\ldots)$ the dimesnsion of
${\bf M}_{ij}(\ldots)$. Observe that ${\bf m}_{ij}(X)\neq 0$ only for $j=i$ and $j=i+1$
and ${\bf m}_{ij}(Y)\neq 0$ only for $j=i+1$.

As ${\bf m}_{0r}(X)=0$, ${\bf m}_{0r}(Y)=0$, for $r\neq 0$ and ${\bf m}_{00}(X)=1$,
${\bf m}_{00}(Y)=1$, it follows ${\bf m}_{00}( \mathcal{L})=0$, i.e. $r < 0$.

The above exact sequence for $j=1$ gives ${\bf M}_{01}( \mathcal{L}) \cong {\bf M}_{11}(X)$
and, as ${\bf m}_{11}(X)=e$, it follows $e=h^0( \mathcal{L}(1))=
h^0(\mathcal{O} _{\projsp{1}}(d+r))$. Then $r+d+1 = e$, i.e. $r=e-d-1$ and $e\le d$.
As $H^0( \mathcal{L})=0$ and $H^0( \mathcal{L}(1))\neq 0$, it follows that $-d <r<0$. i.e. 
$L=\mathcal{T} ^{e-1}(-1)$.
Then the minimal resolution of $\mathcal{L}$ in $\projsp{d}$ is, for $e\neq d$ of the form:
\begin{eqnarray}
0\leftarrow \mathcal{L} \leftarrow e\mathcal{O} (-1) \leftarrow (e-1){d \choose 1}\mathcal{O}(-2) 
\leftarrow (e-2){d \choose 2}\mathcal{O}(-3)\leftarrow \nonumber \\ 
\leftarrow (e-3){d \choose 3}\mathcal{O}(-4) \leftarrow \ldots  
\leftarrow {d \choose e-1} \mathcal{O}(-e) \leftarrow {d \choose e+1}\mathcal{O}(-e-2)
\leftarrow \nonumber \\
\leftarrow 2{d \choose e+2}\mathcal{O} (-e-3) \leftarrow 3{d \choose e+3}\mathcal{O} (-e-4)
\leftarrow \ldots \leftarrow \nonumber \\
\leftarrow (d-e){d \choose d}\mathcal{O}(-d-1) \leftarrow 0 \nonumber
\end{eqnarray}
and the polynomial attached to the minimal resolution of $\mathcal{L}$ in $\projsp{{d+e}}$ is
$(et+(e-1){d \choose 1}t^2x+(e-2){d \choose 2}t^3x^2+\ldots +{d \choose e-1}t^e 
x^{e-1} +
{d \choose e+1}t^{e+2} x^ e + 2{d \choose e+2}t^{e+3}x^{e+1} +\ldots +(d-e-2)
{d \choose d}t^{d+1}
x^{d-1})(1+tx)^e$. This shows that ${\bf m}_{ij}\neq 0$ only for $j=i+1$ if $i<e$ and 
${\bf m}_{ij}\neq 0$ only for $j=i+1$ and $j=i+2$ if $i\ge e$. Then one should have 
the exact 
sequence:
\[
0={\bf M}_{d+e-1,d+e+1}(Y) \leftarrow {\bf M}_{d+e-1,d+e+1}(\mathcal{L}) \leftarrow 
{\bf M}_{d+e,d+e+1}(X)=0 \leftarrow \ldots \ ,
\]
which contradicts ${\bf M}_{d+e-1,d+e+1}(\mathcal{L})\neq 0$. It remains $e=d$ and $r=-1$.

We show now that there are rational curves in $\projsp{{2d}}$, which are  rational normal
curves in a linear subspace $\projsp{d}$ and which have a doubling in $\projsp{{2d}}$ 
with linear syzygies.
For that, consider the Veronese embedding $v$ of $\projsp{2}$ given by the monomials of 
degree
$d$ in the homogeneous coordinates $x_0, x_1, x_2$ on $\projsp{2}$. Then
$u_{i_0,i_1,i_2}=v(x_0,x_1,x_2)$ are the homogeneous coordinates of $\projsp{N}$, $N = 
{d+2 \choose 2}-1$.  Let $P$ be the image of $\projsp{2}$ in $\projsp{N}$.  If we cut  
$P$ with the hyperplane $u_{00d}=0$ we get a {\it l.c.i.} curve $Z$ in $\projsp{{N-1}}$ 
of 
degree $d^2$. The dualizing sheaf
of $Z$ is
\begin{eqnarray}
\omega _Z \cong \omega _P |_ Z \otimes \mathcal{O} _{\projsp{{N-1}}}(1) \quad \hbox{ and 
so }  \nonumber \\
\omega _Z|_{\projsp{1}}\cong \mathcal{O} _{\projsp{1}}(-3+d) \ . \quad \quad \quad 
\nonumber 
\end{eqnarray}

This curve has nilpotents and the reduced structure $X=Z_{red}$ is the rational normal curve 
contained in the linear subspace 
$\mathbb{P} =\projsp{d}$ given by $u_{i_1i_2i_3}=0$, for $i_3\neq 0$ on which the homogeneous 
coordinates are $v_{i_1i_2}=u_{i_1i_20}$ for $i_1+i_2=d$.
Consider the point $u$ ($u_{i_1i_2i_3}=0$, for $(i_1,i_2,i_3)\neq (d00)$). Locally  
in $u$, the curve $Z$ is contained in $N-3$ hypersurfaces which are regular in $x$. 
This shows that $Z$ is a quasiprimitive structure on $X$ (cf. \cite{RefBF}, \cite{RefM1} or
\cite{RefM2}). 
As the degree of $X$ is $d$, the multiplicity of $Z$ in any point is $d$. According to loc. cit. 
the curve $Z$ admits a filtration with Cohen-Macaulay curves $X=Y_1 \subset Y_2 \subset 
\ldots Y_{d-1} \subset Y_d=Z$ which are generically {\it l.c.i.} and such that there are
a line bundle $\mathcal{L}$ on $X$ and $d-2$ effective divisors $D_3$, \ldots $D_d$ on $X$ which
fit into the following exact sequences:
\[
0 \rightarrow \mathcal{L} \rightarrow \mathcal{O} _{Y_2} \rightarrow \mathcal{O} _{Y_1}=
\mathcal{O} _X
\rightarrow 0
\]
\[
0 \rightarrow \mathcal{L}^2(D_2) \rightarrow \mathcal{O} _{Y_3} \rightarrow \mathcal{O} _{Y_2}
\rightarrow 0
\]
\[
0 \rightarrow \mathcal{L}^3(D_2+D_3) \rightarrow \mathcal{O} _{Y_4} \rightarrow \mathcal{O} _{Y_3}
\rightarrow 0
\]
{\hfill \vdots \hfill}
\[
0 \rightarrow \mathcal{L}^{d-1}(D_2+\ldots +D_{d-1} ) \rightarrow \mathcal{O} _{Y_d} \rightarrow 
\mathcal{O} _{Y_{d-1}} \rightarrow 0
\]
Consider $\mathcal{L}|_\projsp{1}=\mathcal{O}_\projsp{1}(r)$. The above exact equences
are possible only for $r\ge -1$.
From the general theory (cf. \cite{RefBF}, \cite{RefM1} or \cite{RefM2}) one has:
\begin{eqnarray}
\omega _Z|_X \cong \omega _X \otimes L^{-(d-1)}(-D_2-\ldots -D_d) \hbox{ and so } 
\nonumber \\
\omega _Z|_{\projsp{1}}\cong \mathcal{O} _{\projsp{1}}(-2-(d-1)r-\delta _2-\ldots - 
\delta _d) \ , \nonumber 
\end{eqnarray}
where $\delta _j\ge 0$ are the degrees of $D_j$'s.

From the two expressions of $\omega _Z|_{\projsp{1}}$ one gets:
\[
-(r+1)(d-1)= \delta _3+\ldots +\delta _d \ .
\]
As the lefthandside member of this equality is $\le 0$, $d \ge 2$ and the right one is 
$\ge 0$, 
it follows:
\[
r=-1 \hbox{ and } \delta _j = 0 \hbox{ for all } j \ .
\]
In other words the multiple structure $Z$ on $X$ is "primitive" in the terminology of 
B\u anic\u a and
Forster. In particular it is l.c.i. and in fact all the $Y_j$'s are l.c.i.\ . $Y_2$ is 
a double structure 
on $X$, of the type we are looking for, because it lies in the linear subspace
$u_{i_1i_2i_3}=0, i_1+i_2+i_3=d,i_3\ge 2$.

In the following we shall find all double curves $Y\subset \projsp{{2d}}$ with support 
$X$ and linear syzygies. To fix the notation, take $X=v_d(\projsp{1})\subset 
\projsp{d}\subset \projsp{{2d}}$ and let $x_0$,$x_1$ be the homogeneous coordinates
on $\projsp{1}$. As already shown, one has to consider doublings $Y$ defined by ideals 
$\mathscr{I} _Y$ which are kernels of surjections $\mathscr{I} _X \rightarrow 
\mathscr{I} _X / \mathscr{I} _X^2 \xrightarrow{p} \mathcal{L}$, where $\mathcal{L}$ is a bundle 
on $X$ such that $\mathcal{L} |_{\projsp{1}}\cong \mathcal{O} _{\projsp{1}}(-1)$. We have immediately
$H^1(\mathscr{I} _Y)=0$.

From the exact cequence:
\[
0\rightarrow H^0(\mathscr{I} _Y(1))\rightarrow H^0(\mathscr{I} _X(1)) 
\rightarrow H^0(\mathcal{L}(1))\rightarrow H^1(\mathscr{I} _Y(1))\rightarrow 0 \ ,
\]
as $\hbox{dim }H^0(\mathscr{I} _X(1))= \hbox{dim }H^0(\mathcal{L}(1))=d$, 
in order to have $Y$ not lying in a hyperplane, one should take $p$ such that
the induced map $H^0(\mathscr{I} _X(1)) \to H^0(\mathcal{L}(1))$ is an isomorphism, i.e.
one should take $p$ such that its component
\[
d\mathcal{O}_{\projsp{{2d}}}(-1)|_X\cong d \mathcal{O}_{\projsp{1}}(-d)\rightarrow
\mathcal{O}_{\projsp{1}}(-1)\cong \mathcal{L}|_{\projsp{1}}
\]
is defined, restricted to $\projsp{1}$, by a basis of the component of degree $d-1$
of the polynomial ring $k[x_0,x_1]$.

From now on we show that such a  $p$ will produce $Y$ with the required properties.
The above exact sequence gives $H^1(\mathscr{I} _Y(1))=0$.
In the commutative diagram:
\[
\begin{CD}
H^0({\mathscr{I}}_X(1))\otimes H^0({\mathcal{O}}_{\projsp{{2d}}}(t))@>>>
H^0(\mathcal{L}(1))\otimes H^0({\mathcal{O}}_{\projsp{{2d}}}(t))\\
@VVV @VVV \\
H^0({\mathscr{I}}_X (t+1)) @>>> H^0(\mathcal{L}(t+1)) 
\end{CD}
\]
our choice of $p$ and the fact that $\mathcal{L}(1)$ is generated by its global sections
show that the bottom morphism is surjective and hence $H^1(\mathscr{I} _Y(t+1))=0$, 
for all $t\ge 0$, i.e. $Y$ is projectively normal.

In the following we shall show that the minimal resolution of $\mathscr{I} _Y$
is pure and all the syzygies are linear. As $Y$ is projectively normal it is
the same to show this property for $\mathcal{O} _Y$.

The minimal resolution of $\mathcal{O} _X$ corresponds to the polynomial 
\begin{eqnarray}
(1+tx)^d(1+{d \choose 2}t^2x+
2{d \choose 3}t^3x^2+\ldots (d-1){d \choose d}t^dx^{d-1})\  = \nonumber \\
(1+tx)^d +1/x(1+tx)^d+dt(1+tx)^{2d-1}-1/x(1+tx)^{2d} \ , \nonumber
\end{eqnarray}
the minimal resolution of $\mathcal{L}$ corresponds to
\begin{eqnarray}
(1+tx)^d(d{d \choose 0}t+(d-1){d \choose 1})t^2x +\ldots +{d \choose d-1}t^dx^{d-1}) \ = 
\nonumber \\
dt(1+tx)^{2d-1} \ .\hspace{6.5cm} \nonumber
\end{eqnarray}
From here it follows that the only nonzero $\bf M_{ij}$ for $X$ and $\mathcal{L}$ are:
\begin{eqnarray}
{\bf M}_{ii}(X) & = &{d\choose i} \quad \hbox{for} \quad 0\le i\le d \nonumber \\
{\bf M}_{i,i+1}(X) & = & {d\choose i+1}+d{2d-1 \choose i}-{2d \choose i+1}\quad 
\hbox{for} 
\quad 1\le i\le 2d-1 \nonumber \\
{\bf M}_{i,i+1}(\mathcal{L})& = & d{2d-1 \choose i} \nonumber \ ,
\end{eqnarray}
where we make the usual convention that the binomial coefficient are zero when they do 
not 
make sense.
As $H^1(\mathcal{L}(\ell))=0$ for any $\ell \ge  0$, from the exact sequence
\[
0 \leftarrow \mathcal{O} _X \leftarrow \mathcal{O} _Y \leftarrow \mathcal{L} \leftarrow 0
\]
one obtains for each $j\ge 0$ the exact sequence (cf. \cite{RefG}):
{\small
\[
0\leftarrow {\bf M}_{0j}(X) \leftarrow {\bf M}_{0j}(Y)\leftarrow {\bf M}_{0j}( \mathcal{L})
\leftarrow {\bf M}_{1j}(X) \leftarrow {\bf M}_{1j}(Y) \leftarrow  {\bf M}_{1j}( \mathcal{L})
\leftarrow \ldots .
\]
}
One gets immediately ${\bf m}_{00}(Y)=1$ and ${\bf m}_{0j}(Y)=0$ for $j\neq 0$ and
from:
\[
0\leftarrow {\bf M}_{01}( \mathcal{L})
\leftarrow {\bf M}_{11}(X) \leftarrow {\bf M}_{11}(Y) \leftarrow  {\bf M}_{11}( \mathcal{L})=0
\]
it follows ${\bf M}_{11}(Y)=0$, i.e. $Y$ is not contained in any hyperplane.
From here  ${\bf M}_{ii}=0$ for any $i\ge 1$.
For $k\ge i+2$ one has the exact sequence:
\[
 0={\bf M}_{ik}(X) \leftarrow {\bf M}_{ik}(Y) \leftarrow  {\bf M}_{ik}( \mathcal{L})=0 \quad ,
\]
which shows that the minimal resolution of $Y$ is linear. As the Hilbert polynomial of 
$Y$ is
that of a rational normal  curve of degree $2d$, this proves the theorem.
\qed

For double rational normal curves one has a similar result to {\bf Theorem (3.c.6)} in
\cite{RefG}. with the notation for the Koszul cohomology from \cite{RefG}:
\begin{theorem}
(Strong Castelnuovo Lemma for double rational normal curves)
\newline (i) Let $\{\Pi _1, \Pi _2,\ldots ,\Pi _e\} \subset 
\projsp{{2d}}$ be a scheme consisting from points and double points (i.e. multiplicity 
2 scheme structure on points), lying on a double rational normal  curve.
Then $\mathcal{K}_{2d-1,1}(\Pi _1, \Pi _2,\ldots ,\mathcal
P_e; \mathcal{O}_{\projsp{{2d}}}(1) ) \ne 0$.
\newline (ii) Conversely, let $\Pi =\{\Pi _1, \Pi _2,\ldots ,\Pi _e \}
\subset \projsp{{2d}}$ be a scheme consisting from points and at most $d$ double points, 
such that, taking the reduced structure one gets points $ P_1, P_2,\ldots , P_e $ 
in linear general position in a linear subspace $\Lambda $ of dimension $d$,
denoted in the following simply $\projsp{}d$. If: 

1) the doubling of $P=\{P_1,\ldots , P_e\}$ takes place outside $\Lambda $, i.e. 
scheme-theoretically $\Lambda \cap \Pi = P$, and 

2) $\mathcal{K}_{2d-1,1}(P) \ne 0$, or equivalently 
$\mathcal{K}_{d-1,1}(P; \projsp{d}) \ne 0$ 

\noindent then there exists a double rational
normal curve which contains the scheme $\Pi $.
\end{theorem}
Proof. (i) The proof is the same as for the corresponding implication in \cite{RefG},
{\bf Theorem (3.c.6)}, nameley if $Y$ is a double rational normal curve containig
$\Pi =\{\Pi _1, \Pi _2,\ldots ,\Pi _e\}$, then the map
\[
\mathcal{K} _{2d-1,1}(Y) \rightarrow \mathcal{K} _{2d-1,1}(\Pi)
\]
is injective. 
\newline (ii) Because of the correspondence between minimal resolutions and Koszul cohomology,
{\bf Lemma 1} shows that $\mathcal{K}_{d-1,1}(P; \projsp{d}) \ne 0$ implies
$\mathcal{K}_{2d-1,1}(P) \ne 0$. Conversely, if 
$\mathcal{K}_{d-1,1}(P; \projsp{d}) = 0$, then ${\bf M}_{d-1,d}(P,\projsp{d})=0$ and so
also ${\bf M}_{d,d+1}(P,\projsp{d})=0$; then ${\bf M}_{2d-1,2d}(P)=0$, i.e. 
$\mathcal{K}_{2d-1,1}(P) = 0$. By this the equivalence in the two conditions in 2) is shown. 

Applying \cite{RefG}, {\bf (3.c.6)}, there is a rational normal curve $X$
in $\projsp{d}$ which contains $P$. It remains to show that there is a doubling $Y$ of $X$,
such that $\Pi \subset Y$. Suppose that the doubled points are $P_1,\ldots ,P_\delta$.
Then $\Pi $ corresponds to an exact 
sequence:
\[
0\rightarrow \mathscr{I} _{\Pi}/\mathscr{I} _P^2 \rightarrow  
I _P/I _P^2 \rightarrow \bigoplus _{\ell =1}^\delta k_\ell \rightarrow 0 \quad ,
\]
where $k_\ell$ is the skyscrapper
$k$ in the point $P_\ell$. We have to show that there is 
a surjection 
\[
N_{X,\projsp{{2d}}}:=\mathscr{I} _X/\mathscr{I} _X^2 \rightarrow \mathcal{L}
\]
such that its restriction to the direct summand $d\mathcal{O} _X (-1)$ of $N_{X,\projsp{{2d}}}$
is also surjective and moreover, the corresponding map $H^0(d\mathcal{O} _X)\to H^0(\mathcal{L}(1))$ is
an isomorphism.
We show this in soundso simple steps.

\noindent {\em Step 1. The map $d\mathcal{O} _X(-1) \rightarrow \mathscr{I} _P/\mathscr{I} _P^2
\rightarrow \bigoplus _{\delta} k$ is surjective.}

\noindent It is enough to show that the map
\[
\mathscr{I} _{\Lambda}/\mathscr{I} _{\Lambda}^2 =d\mathcal{O} _{\Lambda}(-1)\rightarrow
d\mathcal{O} _X(-1) \rightarrow \mathscr{I} _P/\mathscr{I} _P^2 \rightarrow \bigoplus _\delta k
\rightarrow 0
\]
is surjective. This follows from the commutative diagram with exact rows and columns:
\[
\begin{CD}
& & & &{\mathscr{I} _{\Lambda}}/{\mathscr{I} _{\Lambda}^2}
 @= d\mathcal{O} _{\Lambda}(-1)\\
 & & & & @VVV @VVV \\
0 @>>>{\mathscr{I} _{\Pi}}/{\mathscr{I} _P^2} @>>> {\mathscr{I} _P}/{\mathscr{I} _P^2} 
@>>> \bigoplus _{\delta} k@>>> 0\\
& & @|  @VVV @VVV \\
& & {\mathscr{I}_{\Pi}}/{\mathscr{I}_P^2} @>>> 
{\mathscr{I} _P}/{\mathscr{I} _{\Lambda} +\mathscr{I} _P^2} @>>> {\mathscr{I} _P}/{\mathscr{I} 
_{\Pi}+\mathscr{I} _{\Lambda}} & =0
\end{CD}
\]

\noindent {\em Step 2. There exists a map $\mathscr{I} _X/\mathscr{I} _X^2 \to \mathcal{L}$ so that
the induced map $H^0(d\mathcal{O} _X)\to H^0(\mathcal{L}(1))$ is an isomorphism.} 

\noindent We may suppose that
the number of double points in $\Pi$ is exactly $d$, because otherwise we add
new (double) points on $X$. Let $\Delta = \sum _{\ell =1}^d P_\ell$ be the divisor on $X$ 
of points which are the support of the double points in $\Pi$. Apply the functor 
$\mathrm{Hom}(d\mathcal{O} _X(-1), ?)$ to the exact sequence on $X$:
\[
0 \rightarrow \mathcal{L}(-\Delta ) \rightarrow \mathcal{L} \rightarrow \mathcal{O} _\Delta \rightarrow 0 \quad ,
\]
and get the exact sequence:
\begin{eqnarray}
0\rightarrow \mathrm{Hom}(d\mathcal{O} _X(-1), \mathcal{L}(-\Delta))\rightarrow 
\mathrm{Hom}(d\mathcal{O} _X(-1),\mathcal{L})
\rightarrow \mathrm{Hom}(d\mathcal{O} _X(-1), \oplus _d k) \nonumber \\
\rightarrow \mathrm{Ext}^1(d\mathcal{O}_X(-1),\mathcal{L}(-\Delta ))
=dH^1(\mathcal{L}(-\Delta )(1))=dH^1(\mathcal{O} _{\projsp{1}}(-1))=0 \quad \nonumber
\end{eqnarray}
In particular the canonical map $d\mathcal{O} _X(-1) \to \oplus _dk$ factors into
$d\mathcal{O} _X(-1) \buildrel p\over \rightarrow \mathcal{L} \buildrel {evaluation}
\over \longrightarrow \oplus _d k$. Restricted to
$\projsp{1}$, $p$ is given by $d$ homogeneous forms of degree $d-1$. The above composition 
with the evaluation map will produce a $d\times d$ matrix of maximal rank. From here follows 
that the $d$ forms which define $p$ are linearly independent.

\noindent {\em Step 3. The ideal $\mathscr{J} =ker(\mathscr{I} _X \to \mathscr{I} _X/
\mathscr{I} _X^2 \buildrel p\over \rightarrow \mathcal{L})$ defines a double rational normal curve 
with the required properties.}

\noindent Everything follows from the commutative diagram:
\[
\begin{CD}
& & & & 0 @>>>\mathcal{L}(-\Delta ) \\
& & & & @VVV  @VVV \\
0 @>>>  \mathscr{J} @>>> \mathscr{I} _X  @>>> \mathcal{L}  @>>>  0 \\
& &  & &  @VVV @VVV\\
0@>>> \mathscr{I} _{\Pi} @>>> \mathscr{I} _P @>>> \oplus _d k @>>> 0
\end{CD}
\]
with exact rows and columns. \qed

\begin{remark}
The l.c.i. curve $Z$  which we used above has the same Betti numbers as the Veronese 
image of 
$\projsp{2}$. The syzygies of the $Y_j$'s are related by numerous exact sequences of the 
type we used
for $Y_2$ but it seems difficult to find this way the Betti numbers of $v_d(\projsp{2})$.
We can use this idea to produce some particular cases, for instance $d=4$. But, in fact,  
the shape of syzygies for $v_4(\projsp{2})$ can be obtained easier directly:
\end{remark}
\begin{lemma} \label{impure} The minimal resolution of $X:=v_4(\projsp{2}) \hookrightarrow 
\projsp{{14}}$ has the shape:
\begin{eqnarray}
0 \leftarrow \mathcal{O}_X \leftarrow 75\mathcal{O}(-2) \leftarrow 5360\mathcal{O}(-3)
\leftarrow 1947\mathcal{O}(-4) \leftarrow 4488\mathcal{O}(-5) \nonumber \\
\leftarrow 7095\mathcal{O}(-6) 
\leftarrow 7920\mathcal{O}(-7)
\leftarrow 6237\mathcal{O}(-8) \leftarrow 3344\mathcal{O}(-9)\leftarrow 1089\mathcal{O}(-10)
\nonumber \\
\leftarrow 120\mathcal{O}(-11)\oplus 55\mathcal{O}(-12) \leftarrow 24\mathcal{O}(-13)
\leftarrow 3\mathcal{O}(-14)\leftarrow 0 \nonumber
\end{eqnarray}
\end{lemma}
Proof. By \cite{RefBM1} the syzygies have degrees $1$ or $2$, by \cite{RefB} the first 9
matrices have linear entries. Using the Serre duality and the fact that $\omega (1)$ is 
generated by global sections, the minimal resolution should have the shape:
\begin{eqnarray}
0 \leftarrow \mathcal{O}_X \leftarrow 75\mathcal{O}(-2) \leftarrow 5360\mathcal{O}(-3)
\leftarrow 1947\mathcal{O}(-4) \leftarrow 4488\mathcal{O}(-5) \nonumber \\
\leftarrow 7095\mathcal{O}(-6) 
\leftarrow 7920\mathcal{O}(-7)
\leftarrow 6237\mathcal{O}(-8) \leftarrow 3344\mathcal{O}(-9)\leftarrow 1089\mathcal{O}(-10)
\nonumber \\
\leftarrow \begin{matrix}120\mathcal{O}(-11) \\
\oplus \\ (55+x)\mathcal{O}(-12) \end{matrix} \leftarrow  \begin{matrix} x\mathcal{O}(-12) \\
\oplus \\ 24\mathcal{O}(-13) \end{matrix}
\leftarrow 3\mathcal{O}(-14)\leftarrow 0  \quad , \nonumber
\end{eqnarray}
as we can compute step by step. Now use the $SL_3$--invariance of the minimal resolution
of $\omega _X$ to get $x=0$, because the invariant minimal resolution of $\omega _X$ begins 
as follows:
\[
\omega _X \leftarrow V\mathcal {O}(-1) \leftarrow \mathbb{S}_{4,1}(V)\mathcal{O}(-2)
\leftarrow \begin{matrix}(\mathbb{S}_{5,3,1}(V)\oplus \mathbb{S}_{7,1,1}(V))\mathcal{O}(-3)\\
\oplus \\
\mathbb{S}_{8,5}(V)\mathcal{O}(-4)\end{matrix}
\leftarrow \ldots \ ,
\]
as one computes directly. Here $V$ is a vector space of dimmension $3$ and $\mathbb{S}_{\ldots}$
are Schur functors.
\section{Plethysm formulae coming from syzygies}
\label{sec:4}
If we take the action of the group $SL_n({\mathbb C})$ on the projective space in which
$\projsp{n}$ is Veronese there exists an
invariant minimal resolution of the image. This will give recurrence formulae
for $S^m(S^d V)$, where $V$ is such that $\projsp{n}=\mathbb{P}(V^*)$. Our notation is
$\mathbb{P}(V)$ for the space of lines through the origin of $V$. So, when the minimal resolution is
known, the plethysm is
reduced mainly to tensor products, for which one has the \emph{Littlewood-Richardson rule}.

\bigskip
I. We begin with the case of the Veronese embeddings of the line. Let $\projsp{1}=
\mathbb{P}(V^*)$, where $V$ is a vector space of dimension 2.
The Eagon-Northcott complex recalled in {\bf Section }~\ref{sec:2}
can be written invariantly:
\begin{align}
 0 & \leftarrow \mathcal{O} _X \leftarrow \mathcal{O} \leftarrow \Lambda ^2(S^{d-1}(V))\mathcal{O}(-2) 
\leftarrow S^1(V)\cdot \Lambda ^3(S^{d-1}(V))\mathcal{O}(-3) \nonumber \\
& \leftarrow S^2(V)\cdot \Lambda ^4(S^{d-1}(V))\mathcal{O}(-4)\leftarrow S^3(V)\cdot 
\Lambda ^5(S^{d-1}(V))\mathcal{O}(-5) \leftarrow \ldots  \nonumber \\
&\hspace{2cm}  \leftarrow S^{d-3}(V)\cdot \Lambda ^{d-1}(S^{d-1}(V))\mathcal{O}(-d+1)
\nonumber \\
&\hspace{2cm} \leftarrow S^{d-2}(V)\cdot \Lambda ^{d}(S^{d-1}(V))\mathcal{O}(-d)\leftarrow 0 \nonumber
\end{align}
Now, twisting this exact sequence by $\mathcal{O}(\ell )$, with $\ell =2, 3 \ldots $ and
taking the sections one gets recurrence formulae. For the sake of simplicity, we shall denote
$S^t:=S^t(V)$ (then, e.g. $S^r(S^s)\cdot S^p(S^q)$ will replace $S^r(S^s(V))\otimes S^p(S^q(V))$\ )  
and make the convention $S^\ell(\ldots) =0$ for $\ell <0$.
We use also the following formula, valid for $\mathrm{dim}V=2$ (cf. \cite{RefFH}, 
Ex. 11.35 p. 160):
\[
\Lambda ^m(S^m(V))=S^m(S^{n-1+m})
\]
Recall also that the Littlewood-Richardson rule takes for $\mathrm{V}=2$ a very simple shape:
\[
S^mS^n=S^{m+n}+S^{m+n-2}+S^{m+n-4}+\ldots +S^{|m-n|}
\]
Then:

\noindent $\ell =2$:
\[
S^2(S^d)=S^{2d}+S^2(S^{d-2})=0
\]
Applying this repeatedly, one gets (cf. also \cite{RefFH}):
\[
S^2(S^d)=S^{2d}+S^{2d-4}+S^{2d-8}+\ldots 
\]

\noindent $\ell=3$:
\[
S^3(S^d)=S^{3d}+S^2(S^{d-2})\cdot S^d-S^1\cdot S^3(S^{d-3})
\]
\dotfill 

\noindent $\ell=t$:
\begin{align}
& S^t(S^d)= S^{td}+S^2(S^{d-2})\cdot S^{t-2}(S^d)-S^1 \cdot S^3(S^{d-3})\cdot S^{t-3}(S^d) \nonumber \\
& +S^2\cdot S^4(S^{d-4})\cdot S^{t-4}(S^d)-
\ldots +(-1)^{d-1}S^{d-3}\cdot S^{d-1}(S^1)\cdot S^{t-d+1}(S^d) \nonumber \\
& +(-1)^dS^{d-2}\cdot S^d(S^0)\cdot S^{t-d}(S^d) \nonumber \quad .
\end{align}
\begin{example}
\begin{eqnarray}
S^3(S^3) & = & S^9+S^2(S^1)\cdot S^3-S^1\cdot S^3(S^0)=S^9+S^2\cdot S^3-S^1\nonumber \\
& = & S^9+S^5+S^3 \nonumber \\
S^3(S^6) & = & S^{18}+S^2(S^4)\cdot S^6-S^1\cdot S^3(S^3)=S^{18}+(S^8+S^4+I)\cdot S^6 
\nonumber \\
 & & -S^1\cdot (S^9+S^5+S^3)=S^{18}+S^{14}+S^{12}+S^{10}+S^8+2S^6+S^2 \nonumber 
 \end{eqnarray}
etc.
\end{example}
 
\bigskip
II. Consider now the Veronese embedding of $\projsp{2}$ by monomials of degree $2$.
It is a standard exercise to show that the minimal resolution is pure (cf. \cite{RefBM1}, \cite{RefBM2})
and that, invarantly, it has the shape:
\[
0\leftarrow \mathcal{O}_X\leftarrow \mathcal{O}\leftarrow \mathbb{S}_{2,2}\mathcal{O}(-2)
\leftarrow \mathbb{S}_{2,1}\mathcal{O}(-3)\leftarrow \mathbb{S}_{1,1}\mathcal{O}(-4)\leftarrow 0
\]
From here one deduces the following recurrence formula, in the case of $\mathrm{dim}V=3$:
\[
S^t(S^2)=S^{2t}+\mathbb{S}_{2,2}\cdot S^{t-2}(S^2)-\mathbb{S}_{2,1}\cdot S^{t-3}(S^2)+\mathbb{S}_{1,1}
\cdot S^{t-4}(S^2)
\]
\begin{example}
\begin{eqnarray}
S^2(S^2) & = & S^4+\mathbb{S}_{2,2} \nonumber \\
S^3(S^2) & = & S^6+\mathbb{S}_{2,2}\cdot S^1(S^2)-\mathbb{S}_{2,1}\cdot S^0(S^2) = S^6+\mathbb{S}_{4,2}+
\mathbb{S}_{2,2,2} \nonumber \\ 
 & = & S^6+\mathbb{S}_{4,2}+I\nonumber \\
S^4(S^2) & = & S^8+\mathbb{S}_{2,2}\cdot S^2(S^2)-\mathbb{S}_{2,1}\cdot S^2+\mathbb{S}_{1,1}=
\nonumber \\
& & S^8+\mathbb{S}_{6,2}+\mathbb{S}_{4,1}+\mathbb{S}_{4,4}+\mathbb{S}_{3,2}+2\mathbb{S}_2
\nonumber
\end{eqnarray}
\end{example}

\bigskip
III. Take now the Veronese embedding of $\projsp{2}$ by monomials of degree $3$.
The minimal resolution (cf. \cite{RefBM1}, \cite{RefBM2}) can be written invariantly:
\begin{align}
&0\leftarrow \mathcal{O}_X \leftarrow \mathcal{O}\leftarrow \mathbb{S}_{4,2}\mathcal{O}(-2)
\leftarrow (\mathbb{S}_{5,4}+\mathbb{S}_{5,1}+\mathbb{S}_{4,2}+\mathbb{S}_{2,1})\mathcal{O}(-3)
\nonumber \\
&\leftarrow (\mathbb{S}_{6,3}+\mathbb{S}_{5,4}+\mathbb{S}_{5,1}+\mathbb{S}_{4,2}+\mathbb{S}_{3,3}+
S^3+\mathbb{S}_{2,1})\mathcal{O}(-4) \nonumber \\
&\leftarrow (\mathbb{S}_{6,3}+\mathbb{S}_{5,4}+\mathbb{S}_{5,1}+\mathbb{S}_{4,2}+\mathbb{S}_{3,3}+
S^3+\mathbb{S}_{2,1})\mathcal{O}(-5) \nonumber \\
&\leftarrow(\mathbb{S}_{5,4}+\mathbb{S}_{5,1}+\mathbb{S}_{4,2}+\mathbb{S}_{2,1})\mathcal{O}(-6)
\leftarrow \mathbb{S}_{4,2}\mathcal{O}(-7) \leftarrow \mathcal{O}(-9) \leftarrow 0
\nonumber 
\end{align}
From here it follows the following recurrence formula, valid if $\mathrm{dim}V=3$:
\begin{align}
& S^t(S^3)  =   S^{3t}+\mathbb{S}_{4,2}\cdot (S^{t-2}(S^3)-S^{t-7}(S^3))
\nonumber \\
 & -(\mathbb{S}_{5,4}+\mathbb{S}_{5,1}+\mathbb{S}_{4,2}+
\mathbb{S}_{2,1})\cdot (S^{t-3}(S^3)-S^{t-6}(S^3))  \nonumber \\
& +(\mathbb{S}_{6,3}+\mathbb{S}_{5,4}+\mathbb{S}_{5,1}+
\mathbb{S}_{4,2}+\mathbb{S}_{3,3}+
S^3+\mathbb{S}_{2,1})\cdot (S^{t-4}(S^3)-S^{t-5}(S^3)) \nonumber \\
 & +S^{t-9}(S^3) \nonumber 
\end{align}

\bigskip
IV. Consider now the Veronese embedding of $\projsp{2}$ in $\projsp{9}$, by monomials
of degree $3$. Then the minimal resolution (cf. \cite{RefBM1}, \cite{RefBM2}) can be written
invariantely:
\begin{eqnarray}
0\leftarrow \mathcal{O}_X \leftarrow \mathcal{O}\leftarrow \mathbb{S}_{2,2}\mathcal{O}(-2)
\leftarrow \mathbb{S}_{3,2,1}\mathcal{O}(-3)
\leftarrow (\mathbb{S}_{3,3,2}+\mathbb{S}_{3,1})\mathcal{O}(-4) \nonumber \\
\leftarrow \mathbb{S}_{3,2,1}\mathcal{O}(-5) 
\leftarrow\mathbb{S}_{2,2}\mathcal{O}(-6) \leftarrow \mathcal{O}(-8) \leftarrow 0\nonumber 
\end{eqnarray}
and this provides the following recurrence formula, valid for $\mathrm{dim}(V)=4$:
\begin{eqnarray}
S^t(S^2)=S^{2t}+\mathbb{S}_{2,2}\cdot (S^{t-2}(S^2)+S^{t-6}(S^2))-\mathbb{S}_{3,2,1}\cdot 
(S^{t-3}(S^2)+S^{t-5}(S^2))
\nonumber \\
+(\mathbb{S}_{3,3,2}+\mathbb{S}_{3,1})\cdot S^{t-4}(S^2)-S^{t-8}(S^2) \nonumber
\end{eqnarray}
{\bf Acknowledgement.}
A significant part of this paper was written during my visit to Hamburg University
(November 1998-February 1999), so many thanks are due to its hospitality and especially to
Professor Oswald Riemenschneider. Also I want to thank Oldenburg University, where the final
version was written, and especially to Professor Udo Vetter, for the interest in my work and 
for the nice working atmosphere.

\end{document}